\documentclass{amsart}
\usepackage[all]{xy}

\theoremstyle{plain}
\newtheorem{thm}{Theorem}[section]
\newtheorem{cor}[thm]{Corollary}
\newtheorem{lem}[thm]{Lemma}
\newtheorem{prop}[thm]{Proposition}
\newtheorem{fact}[thm]{Fact}
\theoremstyle{definition}

\newtheorem*{thank}{Acknowledgments}

\begin{document}

\title{On the connectedness of Deligne-Lusztig varieties}
\author{Ulrich G\"ortz}
\address{
Mathematisches Institut\\
Beringstr.~1\\
53115 Bonn\\
Germany}
\email{ugoertz@math.uni-bonn.de}
\thanks{The author was partially supported by a Heisenberg grant and by the
SFB/TR 45 ``Periods, Moduli Spaces and Arithmetic of Algebraic Varieties''
of the DFG (German Research Foundation)}

\begin{abstract}
We give a criterion which determines when a union of one-dimen\-sio\-nal
Deligne-Lusztig varieties has a connected closure. We obtain a
new, short proof of the connectedness criterion for Deligne-Lusztig
varieties due to Lusztig.
\end{abstract}

\maketitle

\section{Introduction}

We fix a finite field $\mathbb F_q$ with $q$ elements and of characteristic
$p$, and an algebraic closure $\mathbb F$ of $\mathbb F_q$. Let $G_0$ be a
connected reductive group over $\mathbb F_q$, and let $T_0\subset B_0 \subset
G_0$ be a maximal torus and a Borel subgroup, defined over $\mathbb F_q$.
We denote by $G$, $B$, etc.~the base change to $\mathbb F$, and usually
identify varieties over $\mathbb F$ with their sets of $\mathbb F$-valued
points.
Denote by $W$ the absolute Weyl group, and denote by $\sigma$ the
Frobenius automorphism on $\mathbb F$ (and
$G$, $W$, etc.). Let $S\subset W$ be the set of
simple reflections determined by $B$. For each $w\in W$, we have the
Deligne-Lusztig variety $X(w)$ over $\mathbb F$, defined by
(\cite{deligne-lusztig}, Def.~1.4):
\[
X(w) := \{ g \in G/B;\ g^{-1}\sigma(g) \in
BwB \}.
\]
The unique zero-dimensional Deligne-Lusztig variety is $X({\rm id}) =
(G_0/B_0)(\mathbb F_q)$, and  $X(w)$ is one-dimensional if and only if
$w\in S$. For $s\in S$, the closure $\overline{X(s)}$ of $X(s)$ is
$X(s)\cup X({\rm id})$. See Fact \ref{DL_facts} below.

\begin{thm} \label{thm:connectedness}
Let $I \subseteq S$.  The union 
\[
\overline{X(I)} := X({\rm id}) \cup \bigcup_{s\in I} X(s)
\]
of Deligne-Lusztig varieties is connected if and only if $I$ is not
contained in any proper $\sigma$-stable subset of $S$.
\end{thm}

\begin{proof}
If $I$ is contained in a proper $\sigma$-stable subset of $S$, then
$\overline{X(I)}$ projects onto the set of rational points in the quotient
of $G$ by the corresponding parabolic subgroup, so it cannot be connected.
The converse will be proved in section \ref{sec:proof_connectedness}.
\end{proof}

If $G_0$ is a unitary group, the theorem was proved by Ekedahl and van der
Geer, \cite{ekedahl-vdgeer}, Lemma 7.6 ii), and applied to study the moduli
space of principally polarized abelian varieties over $\mathbb F_q$. See
also \cite{GoertzYu2008-2}, Lemma 7.1 and the discussion following it,
and the remarks below.
For $G_0=GL_n$, the theorem amounts to the statement that one can pass
between any two $\mathbb F_q$-rational flags by a series of steps where in
each step only one subspace of the flag is modified. It is not hard to see
that this is true (and this case is also contained in
\cite{ekedahl-vdgeer}).  

\begin{cor}
For $w\in W$ the following are equivalent:
\begin{enumerate}
\item[(i)] $X(w)$ is irreducible.
\item[(ii)] The closure $\overline{X(w)}$ is connected.
\item[(iii)] The element $w$ is not contained in any $\sigma$-stable proper
standard parabolic subgroup of $W$.
\end{enumerate}
\end{cor}

This result has first been obtained by Lusztig (unpublished); by now, there
are a number of proofs in the literature: by Digne and Michel
\cite{digne-michel} Prop.~8.4, by Bonnaf\'e and Rouquier
\cite{bonnafe-rouquier:irred} Thm.~2, and by Lusztig
\cite{lusztig:perverse-sheaves} Lemma 7.14. The proof given below is
self-contained. For instance, in contrast to \cite{bonnafe-rouquier:irred},
we do not use the fact that Deligne-Lusztig varieties for a Coxeter element
are connected (which was proved by Deligne and Lusztig, see \cite{lusztig})
as an ingredient. Note that it is easy to generalize the criterion to
Deligne-Lusztig varieties in quotients $G/P$ by parabolic subgroups, see
\cite{bonnafe-rouquier:irred}.

\begin{proof}
Since the closure $\overline{X(w)}$ is normal (Fact \ref{DL_facts} (2)),
(i) and (ii) are equivalent.  As in the theorem, it is easy to see that
(ii) implies (iii). Finally assume that (iii) holds, and let $I(w) = \{
s\in S;\ s\le w \}$. The hypothesis says that $I(w)$ is not contained
in a $\sigma$-stable proper subset of $S$, so by the theorem
$\overline{X(I(w))}$ is connected. By Fact \ref{DL_facts} (3), every
connected component of $\overline{X(w)}$ contains a rational point, so
every connected component meets $\overline{X(I(w))}$, and (ii) follows.
\end{proof}

A special case of the theorem is applied by Yu and the author
\cite{GoertzYu2008-2} to study the moduli space (over $\mathbb F_p$) of
$g$-dimensional abelian varieties with Iwahori level structure at $p$. On
this space, one has the Kottwitz-Rapoport stratification, given by the
relative position of the chain of first de Rham cohomology groups and the
chain of Hodge filtrations. This relative position is an element in the
extended affine Weyl group of $GSp_{2g}$. It is striking that
Kottwitz-Rapoport strata share many properties with Deligne-Lusztig
varieties: they are quasi-affine, their closures are normal, the closure
relations are given by the Bruhat order. Moreover, the
proofs are, in a sense, quite similar.
For the quasi-affineness in the case of KR strata, one uses that the moduli
space with Iwahori level structure maps under a finite morphism to a flag
bundle over the moduli space of principally polarized abelian varieties.
On the latter space, the Hodge bundle is ample. The KR strata map to
``relative Deligne-Lusztig varieties'' in this flag bundle. One
shows much in the same way as in Haastert's proof of the
quasi-affineness of Deligne-Lusztig varieties (see \cite{haastert:diplom})
that the strata are quasi-affine; see \cite{GoertzYu2008-2}, Thm.~5.4.
To determine the closure relations, one can use a ``local model diagram'': in
the case of Deligne-Lusztig varieties, the local structure is the same as
the local structure of Schubert cells (and for their closures etc., see
the proof of Fact \ref{DL_facts} (1), (2) below). For Kottwitz-Rapoport
strata, the local structure is described by the local model (see Rapoport
and Zink's book \cite{rapoport-zink}, Ch.~3) which can be identified with a
union of Schubert varieties in an affine flag variety \cite{goertz:sympl}.
See also the survey articles by Rapoport \cite{rapoport:guide} and Haines
\cite{haines:clay}.

Those Kottwitz-Rapoport strata which are contained in the supersingular
locus are actually disjoint unions of Deligne-Lusztig varieties (by
\cite{GoertzYu2008} Cor.~6.5 and \cite{GoertzYu2008-2} Cor.~7.4). Using the
theorem above, one can show that all KR strata which are not contained in
the supersingular locus, are connected (and conversely, if the level
structure away from $p$ is small enough), \cite{GoertzYu2008-2} Thm.~7.3.
See also \cite{ekedahl-vdgeer} for a similar application. The reasoning in
the proof of the Corollary also appears in Oort's paper \cite{oort} about
the canonical stratification (now called the Ekedahl-Oort stratification)
of the moduli space of principally polarized abelian varieties.

One can hope that this approach also yields results about the connectedness
of closed affine Deligne-Lusztig varieties in the Iwahori case (see
Viehmann \cite{viehmann} for the case of a maximal parahoric subgroup).
However, an additional effort is required because very little is known about
properties of affine Deligne-Lusztig varieties, for instance, whether each
connected component contains a rational point in its closure.

\begin{thank}
I am grateful to Tetsushi Ito for his comments on the text, and to George
Lusztig for pointing me to the reference \cite{lusztig:perverse-sheaves}.
\end{thank}

\section{Basic properties of Deligne-Lusztig varieties}

We regard the set $(G_0/B_0)(\mathbb F_q)$ as a subset of $G/B$; it is the
set of points fixed by $\sigma$, and as usual we call its elements the
rational points of $G/B$.

On $W$, we have the Bruhat order $\le$, and the length function $\ell$.  We
denote by $C_v := BvB/B$ the Schubert cell associated with $v\in W$. Its
closure $\overline{C}_v$ is called the Schubert variety for $v$.  The set
$C_v \cap (G_0/B_0)(\mathbb F_q)$ is non-empty if and only if $v$ is fixed by
$\sigma$. In this case $C_v$ is defined over $\mathbb F_q$, and is
isomorphic to $\mathbb A^{\ell(v)}$ over $\mathbb F_q$, and we denote its set of rational
points by $C_v(\mathbb F_q)$.
We denote the subgroup of $W$ of elements fixed by $\sigma$ by $W^\sigma$.

We need the following well-known properties of Deligne-Lusztig varieties:

\begin{fact} \label{DL_facts}
\begin{enumerate}
\item
Let $w\in W$.
We have $\overline{X(w)} = \bigcup_{v\le w} X(v)$.
\item
Let $v\le w$ be elements of $W$. Fix points $v_1 \in X(v)(\mathbb F)$,
$v_2\in C_v(\mathbb F)$. Then the singularities of $\overline{X(w)}$ in
$v_1$ and of $\overline{C}_w$ in $v_2$ are smoothly equivalent.
In particular, $\overline{X(w)}$ is normal.
\item
Let $w\in W$, and
let $X \subseteq X(w)$ be an irreducible component. Then the closure
$\overline{X}$ contains a rational point of $G/B$.
\end{enumerate}
\end{fact}

\begin{proof}
To make this paper self-contained, we sketch proofs of these facts.
For (1), (2), we repeat the proof given in \cite{GoertzYu2008}, Section 5.
Let $L\colon G\rightarrow G$ be the Lang map $g\mapsto g^{-1}\sigma(g)$.
The composition $\xymatrix{G \ar[r]^L & G \ar[r]^{{\rm pr}} & G/B}$ is
$B$-equivariant if $B$ acts in $G$ by right multiplication, and on $G/B$ by
$b.gB = b^{-1}gB$. Taking the quotient, we get a smooth morphism of
relative dimension $\dim B$
\[
G/B \rightarrow [B\backslash G/B]
\]
from $G/B$ to the stack quotient $[B\backslash G/B]$. The underlying
topological space of this stack is just the Weyl group $W$, and the fiber
over $w\in W$ is $X(w)$. This proves the first two assertions.

To prove (3), it is enough to show (because of (1)) that $G_0(\mathbb F_q)$
acts transitively on the set of connected components of $X(w)$. We recall
the argument presented in \cite{bonnafe-rouquier:irred}: the Lang map $L$
is a $G_0(\mathbb F_q)$-torsor, so $G_0(\mathbb F_q)\backslash L^{-1}(BwB)
\cong BwB$ is connected. On the other hand, $G_0(\mathbb F_q)\backslash
L^{-1}(BwB)$ surjects onto $G_0(\mathbb F_q)\backslash X(w)$.
Another way to prove (3) is to use that all Deligne-Lusztig varieties are
quasi-affine, as was proved by Haastert \cite{haastert:diplom}. This
makes the analogy with \cite{GoertzYu2008-2} even closer.
\end{proof}

\section{Proof of the connectedness criterion}
\label{sec:proof_connectedness}

\begin{lem} \label{key-lemma}
Let $s\in S$.
\begin{enumerate}
\item
Let $v\in W$. Then $X(s) \cap C_v$ is empty, unless $v \in W^\sigma$ and
$vs<v$, in which case the intersection is equi-dimensional of dimension $1$.
\item
Let $v\in W^\sigma$, such that $vs<v$, and let $x \in C_v(\mathbb F_q)$. 
Then there exists an irreducible component $X \subset X(s)$ which is
contained in $C_v$ and whose closure $\overline{X}$ (which is an
irreducible component of $\overline{X(s)}$) contains $x$.
Furthermore, $\overline{X}$ contains a rational point of some cell
$C_{v'}$, $v'\in W^\sigma$, $v'<v$.
\end{enumerate}
\end{lem}

\begin{proof}
Denote by $U$ the unipotent radical of $B$, and by $U^-$
the unipotent radical of the Borel opposite to $B$. We let $U_v =
vU^-v^{-1}\cap U$, and then have an isomorphism $U_v\rightarrow C_v$,
$g\mapsto gv$. It induces an isomorphism
\[
\xymatrix{
    X(s) \cap C_v \ar[r]^<<<<{\cong} & 
    \{ g \in U_v;\ (gv)^{-1}\sigma(gv) \in BsB \} =: U_v(s)
}
\]
The Lang map $L\colon u\mapsto u^{-1}\sigma(u)$ for $U_v$ induces a finite
\'etale morphism
\[
\xymatrix{
    U_v(s) \ar[r]^<<<<L &
vBsB\sigma(v^{-1}) \cap U_v \ar[r]^<<<<{\cong} &
    BsB\sigma(v^{-1})v \cap v^{-1}Uv \cap U^-.
}
\]
The second morphism is conjugation by $v^{-1}$.
If $\sigma(v) \ne v$, then $v^{-1}\not\in BsB\sigma(v^{-1})B$, and it
follows that in this case the intersection $BsB\sigma(v^{-1})v \cap
v^{-1}Uv$ and a fortiori the right hand side is empty. 

Now assume that $v\in W^\sigma$. The Bruhat decomposition shows that the
right hand side is non-empty if and only if $vs<v$, and that it is
irreducible of dimension $1$ in the latter case. Because $L$ is finite
\'etale, (1) is proved.

Now consider $v$, $x$ as in (2).
Let $X'\subset X(s)$ be an irreducible component which meets $C_v$. It
follows from part (1) that $X' \cap C_v$ has dimension $1$. We claim that
$X' \subseteq C_v$. Otherwise, the zero-dimensional complement $X'
\setminus C_v$ would intersect some other Schubert cell $C_{v'}$, but this is
impossible by (1), applied to $v'$.  Furthermore, denoting by
$\overline{X'}$ the closure of $X'$, we have, similarly as above, a finite
(\'etale) map
\[
\overline{X'} \cap C_v \subset \overline{X(s)} \cap C_v \longrightarrow
(BsB\cup B) \cap v^{-1}Uv \cap U^- = (BsB \cap v^{-1}Uv \cap U^-) \cup \{
1 \}.
\]
Its image is one-dimensional and closed, and hence equal to the right hand
side. In particular, $\overline{X'}\cap C_v(\mathbb F_q)$, which is the
fiber over $1$, is non-empty. By applying a suitable element of $B(\mathbb
F_q)$ to $X'$ we produce an irreducible component $X \subset X(s)$ which is
contained in $C_v$ and whose closure contains $x$. 

Because $\overline{X}$ is projective, it is not contained in $C_v$, so it
intersects non-trivially with some cell $C_{v'}$ in the closure of $C_v$;
since $X \subset C_v$, all points in $\overline{X}\setminus C_v$ lie in
$X({\rm id})$, i.~e.~are rational points of $G/B$.
\end{proof}

We note that part (1) of the lemma, except for the equi-dimensionality
statement, follows from (a very special case of) Prop.~8.2 in
\cite{digne-michel}, and the proof given above is basically the relevant
part of the proof in loc.~cit.

\begin{prop} \label{prop}
Let $I\subseteq S$ be a subset which is not contained in a $\sigma$-stable
proper subset of $S$.  Let $x \in (G_0/B_0)(\mathbb F_q)$ be a rational point,
say $x\in C_v$. Assume $v\ne {\rm id}$.  Then there exists $s\in I$ such
that $vs < v$. For each such $s$, there exist $v'<v$ and a rational point
in $C_{v'}(\mathbb F_q)$ which lies in the same connected component of
$\overline{X(I)}$ as $x$.
\end{prop}

\begin{proof}
To show the existence of $s$, note that $C_v\cap (G_0/B_0)(\mathbb F_q)\ne
\emptyset$ implies that $v\in W^\sigma$, so in particular, $vs< v$ if and
only if $v\sigma(s)< v$. This implies that unless $v={\rm id}$, there
exists $s\in I$ as desired.
Now we apply part (2) of the lemma.
\end{proof}

By induction, the proposition implies that every rational point of
$(G_0/B_0)(\mathbb F_q)$ is in the same connected component of $\overline{X(I)}$ as
the point $C_{\rm id}$. Because  every irreducible component of
$\overline{X(I)}$ contains a rational point (Fact \ref{DL_facts} (3)), the
theorem is proved.

\section{More detailed analysis}

Let $s\in S$, and let $X_1 \subseteq X(s)$ be the irreducible component
whose closure contains $C_{\rm id}$. If the group $G$ is split, the
situation is relatively simple: $\overline{X_1} = \overline{C}_s$. In
particular, all irreducible components of $\overline{X(s)}$ are just
projective lines, and contain $q+1$ rational points. The finite group
$G_0(\mathbb F_q)$ acts on $\overline{X(s)}$, and given $v$ with $vs < v$,
$vsX_1$ is an irreducible component of $\overline{X(s)}$ which connects a
rational point in $C_v$ with one in $C_{vs}$, as in the proposition. We
also see that $X(s)$ has $\#(G_0/B_0)(\mathbb F_q)/(1+q)$ irreducible
components. In general, the picture is to some extent similar:

For $s\in S$, we denote by $W^s$ the subgroup generated by all elements
$\sigma^i(s)$, $i\in \mathbb Z$, and by $w_0^s$ its longest element. In
particular, $w_0^s = s$ if $G$ is split. In every case, $w_0^s\in
W^\sigma$. We need to know the structure of $W^\sigma$, as described in

\begin{lem} {\rm (Steinberg \cite{steinberg} \S 1; Kottwitz, Rapoport
\cite{kottwitz-rapoport:alcoves} Prop.~2.3)}
The subgroup $W^\sigma\subseteq W$ together with the set $\{ w_0^s;\ s\in
S\}$ is a Coxeter system. The Bruhat order on $W^\sigma$ is the restriction
of the Bruhat order on $W$ to $W^\sigma$.
\end{lem}

In particular, the lemma implies that $W^s \cap W^\sigma = \{ {\rm id},
w_0^s \}$, and that for $v\in W^\sigma$, $vs<v$ if and only if $vw_0^s <
v$.  Now it is easy to extend the description given above for split groups
to the general case.

\begin{prop}
Let $s\in S$.
\begin{enumerate}
\item
Let $X_1 \subset X(s)$ be the irreducible component whose closure contains
the point $C_{\rm id}$. Then $\overline{X_1} \cap (G_0/B_0)(\mathbb F_q) =
C_{\rm id} \cup C_{w_0^s}(\mathbb F_q)$.
\item
Let $v \in W^\sigma$ such that $vs < v$. Then $vw_0^s < v$. For every
$x\in C_v(\mathbb F_q)$, there exists an irreducible component of $X(s)$
whose closure contains $x$ and a rational point of $C_{vw_0^s}$.
\end{enumerate}
\end{prop}

Furthermore, we get a formula for the number of connected components. We
extend the notation introduced above as follows: For $w\in W$, the subgroup
$W^w$ denotes the smallest $\sigma$-stable standard parabolic subgroup of
$W$ which contains $w$, i.~e.~the subgroup generated by all $W^s$, $s\le
w$. We denote by $P^w_0$ the standard parabolic subgroup of $G_0$
corresponding to $W^w$.  For a subset $H \subseteq W$, let $N(H) = \sum_{w
\in H \cap W^\sigma} q^{\ell(w)}$. For instance, $N(W) =
\#(G_0/B_0)(\mathbb F_q)$, and more generally $N(W^w) =
\#(P^w_0/B_0)(\mathbb F_q)$. Concretely, for $s\in S$ we have $N(W^s) = 1
+ q^{\ell(w_0^s)}$.

\begin{prop}
Let $w\in W$. The projection $G/B \rightarrow G/P^w$ restricts to a
surjection $X(w) \rightarrow (G_0/P^w_0)(\mathbb F_q)$, whose fibers are
the connected components of $X(w)$. In particular, $X(w)$ has $N(W)/N(W^w)$
connected components, and every connected component has $N(W^w)$ rational
points of $G/B$ in its closure.
\end{prop}

\begin{proof}
This is an easy consequence of the above. See \cite{bonnafe-rouquier:irred}
or \cite{GoertzYu2008}, Cor.~5.3.
\end{proof}


\end{document}